\documentclass[12pt, twoside]{article}
\usepackage{amsmath,amsthm,amssymb}
\usepackage{times}
\usepackage{enumerate}
\usepackage{amssymb}
\usepackage{graphicx}

\def\titlerunning#1{\gdef\titrun{#1}}
\makeatletter
\def\author#1{\gdef\autrun{\def\and{\unskip, }#1}\gdef\@author{#1}}

\makeatother

\newtheorem{thm}{Theorem}[section]
\newtheorem{cor}[thm]{Corollary}
\newtheorem{lem}[thm]{Lemma}

\theoremstyle{definition}

\newtheorem{rem}[thm]{Remark}

\numberwithin{equation}{section}

\begin{document}

\baselineskip=17pt

\titlerunning{Existence of bounded Palais-Smale sequences}

\title{Some abstract results on the existence of bounded Palais-Smale sequences}

\author{
Michela Guida,
\ \ Sergio Rolando}

\date{
\begin{footnotesize}
\emph{
Dipartimento di Matematica ``Giuseppe Peano'' \smallskip\\
Universit\`{a} degli Studi di Torino, Via Carlo Alberto 10, 10123 Torino, Italy \\
e-mail:} michela.guida@unito.it, sergio.rolando@unito.it
\end{footnotesize}
}
\maketitle

\begin{abstract}
Without compactness assumptions, we prove some abstract results which show
that a $C^{1}$ functional $I:X\rightarrow \mathbb{R}$ on a Banach space $X$
admits bounded Palais-Smale sequences provided that it exhibits some
geometric structure of minimax type and a suitable behaviour with respect to
some sequence of continuous mappings $\psi _{n}:X\rightarrow X$. This work
is a preliminary version of a forthcoming paper, where applications to
nonlinear equations without Ambrosetti-Rabinowitz type assumptions will also
be given.
\end{abstract}

\section{Introduction}

This paper is concerned with the existence of bounded Palais-Smale sequences
for $C^{1}$ functionals defined on a Banach space and satisfying suitable
minimax conditions. This kind of study is a fundamental step in minimax
methods of critical point theory, where existence of critical points for a
given functional $I:X\rightarrow \mathbb{R}$ on a Banach space $X$ is usually
obtained by searching for sequences of ``almost critical'' points and then
showing their convergence to an exact critical point. More precisely, a
minimax method commonly conforms to the following scheme (see \cite{Jabri}, 
\cite{Willem}):

\begin{itemize}
\item  a geometric condition is required, involving a relation between the
values of $I$ over sets that satisfy a topological intersection property (%
\textsl{linking property});

\item  by a quantitative deformation lemma (or by Ekeland's variational
principle), the existence of a \textsl{Palais-Smale sequence} for $I$ is
provided, namely, a sequence $\left\{ w_{n}\right\} \subset X$ such that 
\[
I\left( w_{n}\right) \rightarrow c\quad \text{and}\quad I^{\prime }\left(
w_{n}\right) \rightarrow 0~\text{in }X^{\prime }
\]
where $c$ is some real value defined by a minimax procedure (and $X^{\prime }
$ is the dual space of $X$);

\item  some \textit{a posteriori} compactness property of $\left\{
w_{n}\right\} $ is proved, with a view to concluding that the ``almost
critical'' points $w_{n}$ actually approximate an exact critical point of $I$.
\end{itemize}

\noindent A wide range of abstract results describing geometric structures
of $I$ for which the first two steps of the above scheme can be successfully
carried out is well known in the literature and the first and most popular
one certainly concerns the mountain-pass geometry, introduced by Ambrosetti
and Rabinowitz in their renowed paper \cite{Ambr-Rab}. Then, under the
assumption that all the Palais-Smale sequences of $I$ admit a strongly
convergent subsequence (\textsl{Palais-Smale condition}), the minimax level $%
c$ is a critical value for $I$ (assuming $I$ of class $C^{1}$).
Unfortunately, the Palais-Smale condition fails in many concrete situations
and one has to carry out the final step of the method directly, by showing
that a particular Palais-Smale sequence $\left\{ w_{n}\right\} $ can be
found, which is relatively compact. A first and key point in this direction
is commonly to establish the boundedness of $\left\{ w_{n}\right\} $, which,
if $X$ is reflexive, immediately yields that (up to a subsequence) it weakly
converges to some $w\in X$. Then $w$ is a critical point if one succeeds in
bringing the amount of compactness needed to conclude strong convergence, or
even just in showing that $I^{\prime }\left( \cdot \right) \varphi $ is
weakly continuous for every fixed $\varphi \in X$ (to be precise, some more
work is often required, since one usually needs $w\neq 0$).

The problem of obtaining bounded Palais-Smale sequences has been faced by
several authors in different specific contexts (see, for instance, \cite
{Miya-Souto}, \cite{Schech-Zou 04}, \cite{Wang-Tang}, \cite{Zhou}, \cite
{Zhou 98} and the references contained in \cite{JJ}) and it has been
sistematically approached by Jeanjean in \cite{JJ} (see also \cite{JJ 2}-%
\cite{JJ-Toland}), where a general method deriving from Struwe's \textsl{%
monotonicity trick} \cite{Strick} is presented. More precisely, the author
formulates an abstract result establishing very general conditions in order
that, given a family of functionals depending on a real parameter and
exhibiting a ``uniform'' mountain-pass geometry, almost every functional of
the family has a bounded Palais-Smale sequence at its mountain-pass level.
Then, roughly speaking, the method consists in exploiting such a result
(together with some additional compactness properties) in order to obtain a
special sequence of ``almost critical'' points for a given functional,
consisting of exact critical points of nearby functionals and thus
possessing extra properties which can help in proving its boundedness. Other
results in a similar abstract spirit can be found in \cite{Ambr-Ruiz}, \cite
{Rabier}, \cite{Schech-Zou}, \cite{Szulkin-Zou}, \cite{Will-Zou}, \cite{Zou}.

Here we present some abstract results which show that, for different
geometric structures of the functional $I:X\rightarrow \mathbb{R}$ (and without
any compactness assumption), the second step of the above scheme can be
performed in such a way that a bounded Palais-Smale sequence is directly
yielded, provided that $I$ exhibits a suitable behaviour with respect to
some sequence of continuous mappings $\psi _{n}:X\rightarrow X$.

Precisely, we consider a Banach space $\left( X,\left\| \cdot \right\|
\right) $, a functional $I\in C^{1}\left( X,\mathbb{R}\right) $ of the form 
\begin{equation}
I\left( u\right) =A\left( u\right) -B\left( u\right) \text{\quad for all }%
u\in X  \label{I=A-B}
\end{equation}
and we assume that:

\begin{itemize}
\item[\textbf{(}$\mathbf{A}$\textbf{)}]  $A$ is nonnegative and such that 
\[
\forall l_{2}>l_{1}>0,\quad \exists \delta _{l_{1},l_{2}}>0,\quad \forall
u\in X,\quad \mathrm{dist}\left( u,A^{l_{1}}\right) \leq \delta
_{l_{1},l_{2}}~\Rightarrow ~u\in A^{l_{2}}
\]
(where $A^{l}:=\left\{ u\in X:A\left( u\right) \leq l\right\} $ and $\mathrm{%
dist}\left( u,A^{l_{1}}\right) :=\inf_{h\in A^{l_{1}}}\left\| u-h\right\| $);

\item[\textbf{(}$\mathbf{\Psi }$\textbf{)}]  there exists a sequence of
mappings $\left\{ \psi _{n}\right\} \subset C\left( X,X\right) $ such that $%
\forall n$ there exist $\alpha _{n}>\beta _{n}>0$ satisfying

\begin{itemize}
\item[\textbf{(}$\mathbf{\Psi }$\textbf{{\tiny 1}) }]  $A\left( u\right)
\geq \alpha _{n}A\left( \psi _{n}\left( u\right) \right) $ and $B\left(
u\right) \leq \beta _{n}B\left( \psi _{n}\left( u\right) \right) $ for all $%
u\in X$\smallskip 

\item[\textbf{(}$\mathbf{\Psi }$\textbf{{\tiny 2}) }]  $\lim\limits_{n%
\rightarrow \infty }\alpha _{n}=\lim\limits_{n\rightarrow \infty }\beta
_{n}=1$ and $\liminf\limits_{n\rightarrow \infty }\dfrac{\left| 1-\beta
_{n}\right| }{\alpha _{n}-\beta _{n}}<\infty .$
\end{itemize}
\end{itemize}

\noindent Observe that \textbf{(}$\mathbf{A}$\textbf{)} is quite a mild
assumption in practice, since in many applications one has $A\left( u\right)
=\left( \mathrm{const.}\right) \left\| u\right\| ^{p}$, $p>0$, so that for
every $u,v\in X$ it holds 
\[
A\left( u\right) \leq \left( \mathrm{const.}\right) \left( \left\| v\right\|
+\left\| u-v\right\| \right) ^{p}\rightarrow A\left( v\right) \quad \text{as 
}\left\| u-v\right\| \rightarrow 0. 
\]
Notice that in this case the existence of a bounded Palais-Smale sequence is
obviously yielded by the existence of a Palais-Smale sequence on which $A$
is bounded, but such an implication holds also true if $A$ or $\left|
B\right| $ is coercive with respect to the norm of $X$. Furthermore, we
point out that, even when neither $A$ nor $\left| B\right| $ is coercive
with respect to the norm, the existence of a Palais-Smale sequence on which $%
A$ is bounded can be a relevant information in order to get a Palais-Smale
sequence which is bounded in $X$ (see \cite{BBR 3}).

Under the above assumptions \textbf{(}$\mathbf{A}$\textbf{)} and \textbf{(}$%
\mathbf{\Psi }$\textbf{)}, we will prove that the presence of some geometric
structure of minimax type is essentially sufficient in order that $I$
exhibits, at the minimax level $c$, a Palais-Smale sequence on which $A$ is
bounded, i.e., a sequence $\left\{ w_{n}\right\} \subset X$ satisfying 
\begin{equation}
I\left( w_{n}\right) \rightarrow c\,,\quad I^{\prime }\left( w_{n}\right)
\rightarrow 0~\text{in }X^{\prime }\,,\quad \sup_{n}A\left( w_{n}\right)
<\infty \,.  \label{BDD PS}
\end{equation}
This is for instance the case when $I$ has a mountain-pass geometry, as the
following theorem says.

\begin{thm}
\label{THM: MP}Let $\left( X,\left\| \cdot \right\| \right) $ be a Banach
space and $I\in C^{1}\left( X,\mathbb{R}\right) $ a functional of the form (\ref
{I=A-B}) satisfying \textbf{\emph{(}}$\mathbf{A}$\textbf{\emph{)}} and 
\textbf{\emph{(}}$\mathbf{\Psi }$\textbf{\emph{)}}. If there exist $r>0$ and 
$\bar{u}\in X$ with $\left\| \bar{u}\right\| >r$ such that 
\begin{equation}
\inf_{\left\| u\right\| =r}I\left( u\right) >I\left( 0\right) \geq I\left( 
\bar{u}\right) \,~\text{and}~~\lim\limits_{n\rightarrow \infty }\left\| \psi
_{n}\left( 0\right) \right\| =\lim\limits_{n\rightarrow \infty }\left\| \psi
_{n}\left( \bar{u}\right) -\bar{u}\right\| =0  \label{Thm(MP): hp}
\end{equation}
then there exists a sequence $\left\{ w_{n}\right\} \subset X$ satisfying (%
\ref{BDD PS}), where 
\begin{equation}
c:=\inf\limits_{\gamma \in \Gamma }\max\limits_{u\in \gamma \left( \left[
0,1\right] \right) }I\left( u\right) ,~~\Gamma :=\left\{ \gamma \in C\left(
\left[ 0,1\right] ,X\right) :\gamma \left( 0\right) =0,\,\gamma \left(
1\right) =\bar{u}\right\} .  \label{Thm(MP): def}
\end{equation}
\end{thm}

Theorem \ref{THM: MP} has been already announced and used in \cite{BR TMA}
and it will be obtained in Section \ref{SEC: abstract} as a consequence of a
more general minimax principle (Theorem \ref{THM MAIN}), which also allows
to deduce similar results concerning geometries of saddle-point or linking
type (Theorems \ref{THM: SD} and \ref{THM: L}). The arguments leading to
such a principle derive from the ones of \cite{Ambr-Str}, do not employ the
monotonicity trick and have been previously used in \cite{BBR 3} and \cite
{GR curves} within specific contexts. An estimate of the bound of $A\left(
w_{n}\right) $ can also be obtained (cf. Remark \ref{Rmk: bound}).

This work is a preliminary version of a forthcoming paper \cite{GR bdd},
where applications to nonlinear equations without Ambrosetti-Rabinowitz type
assumptions will also be given.

\section{Abstract results \label{SEC: abstract}}

In this section we give our abstract results concerning the existence of
Palais-Smale sequences satisfying (\ref{BDD PS}) in the functional framework
described in the Introduction. Accordingly, throughout the section we assume
that $\left( X,\left\| \cdot \right\| \right) $ is a Banach space and $I\in
C^{1}\left( X,\mathbb{R}\right) $ is a functional of the form (\ref{I=A-B})
satisfying \textbf{(}$\mathbf{A}$\textbf{)} and \textbf{(}$\mathbf{\Psi }$%
\textbf{)}.

Our main result is the following.

\begin{thm}
\label{THM MAIN}Let $\left( M,d\right) $ be a metric space and $M_{0}\subset
M$ a compact subspace such that there exist $\varepsilon _{0}>0$ and $\sigma
\in C\left( M,M\right) $ satisfying

\begin{itemize}
\item[i) ]  $\sigma $ is uniformly continuous on $M_{0}^{\varepsilon
_{0}}:=\left\{ p\in M:d\left( p,M_{0}\right) \leq \varepsilon _{0}\right\} $

\item[ii) ]  $\sigma \left( M_{0}^{\varepsilon _{0}}\right) \subseteq M_{0}$

\item[iii) ]  $\sigma _{\mid _{M_{0}}}=\mathrm{id}$ .
\end{itemize}

\noindent Let $\Gamma _{0}\subset C\left( M_{0},X\right) $ be such that $%
\bigcup_{\gamma _{0}\in \Gamma _{0}}\gamma _{0}\left( M_{0}\right) $ is
compact and define $\Gamma :=\{\gamma \in C\left( M,X\right) :\gamma _{\mid
_{M_{0}}}\in \Gamma _{0}\}$. Assume that 
\begin{equation}
\sup_{\gamma _{0}\in \Gamma _{0}}\max_{u_{0}\in \gamma \left( M_{0}\right)
}I\left( u\right) <c:=\inf_{\gamma \in \Gamma }\sup_{u\in \gamma \left(
M\right) }I\left( u\right) <\infty   \label{a < c}
\end{equation}
and 
\begin{equation}
\lim\limits_{n\rightarrow \infty }\sup_{\gamma _{0}\in \Gamma
_{0}}\max_{u\in \gamma _{0}\left( M_{0}\right) }\left\| \psi _{n}\left(
u\right) -u\right\| =0\,.  \label{sup max ->0}
\end{equation}
Then there exists a sequence $\left\{ w_{n}\right\} \subset X$ satisfying (%
\ref{BDD PS}).
\end{thm}

\begin{rem}
\label{Rmk: bound}According to the proof of Theorem \ref{THM MAIN}, what we
obtain exactly is that $\forall a>b:=\left| c\right| \liminf_{n\rightarrow
\infty }\left| 1-\beta _{n}\right| /\left( \alpha _{n}-\beta _{n}\right) $
the functional $I$ has a Palais-Smale sequence $\left\{ w_{n}\right\} $ at
level $c$ such that $A\left( w_{n}\right) \leq a$ for all $n$. Observe that
it is thus likely to occur that $A\left( w_{n}\right) \rightarrow 0$ if $b=0$%
.
\end{rem}

Since quite long and technical, we displace the proof of Theorem \ref{THM
MAIN} in the next section and first derive some consequences of it,
beginning with the proof of Theorem \ref{THM: MP}.\bigskip

\proof[Proof of Theorem \ref{THM: MP}]
Set $M:=\left[0,1\right] \subset \mathbb{R}$, $M_{0}:=\left\{ 0,1\right\} $, $\Gamma
_{0}:=\left\{ \gamma _{0}\right\} $ where $\gamma _{0}\in C\left(
M_{0},X\right) $ is defined by $\gamma _{0}\left( 0\right) :=0$ and $\gamma
_{0}\left( 1\right) :=\bar{u}$, and 
\[
\sigma \left( p\right) :=\left\{ 
\begin{array}{ll}
0\medskip & \text{if }p\in \left[ 0,\varepsilon _{0}\right] \\ 
\dfrac{p-\varepsilon _{0}}{1-2\varepsilon _{0}}\medskip & \text{if }p\in
\left[ \varepsilon _{0},1-\varepsilon _{0}\right] \\ 
1 & \text{if }p\in \left[ 1-\varepsilon _{0},1\right]
\end{array}
\right. 
\]
for any fixed $\varepsilon _{0}\in \left( 0,1/2\right) $. So $\Gamma $ and $%
c $ of Theorem \ref{THM MAIN} coincide with the ones of (\ref{Thm(MP): def})
and (\ref{a < c})-(\ref{sup max ->0}) are fulfilled by (\ref{Thm(MP): hp}).
The conclusion thus follows by applying Theorem \ref{THM MAIN}.%
\endproof%

\begin{thm}
\label{THM: SD}Assume $X=Y\oplus Z$ with $\dim Y<\infty $ and assume that
there exists $\rho >0$ such that 
\begin{equation}
\inf_{u\in Z}I\left( u\right) >\max_{u\in M_{0}}I\left( u\right) \quad \text{%
and}\quad \lim\limits_{n\rightarrow \infty }\max_{u\in M_{0}}\left\| \psi
_{n}\left( u\right) -u\right\| =0  \label{Thm(SD): hp}
\end{equation}
where $M_{0}:=\left\{ u\in Y:\left\| u\right\| =\rho \right\} $. Then there
exists a sequence $\left\{ w_{n}\right\} \subset X$ satisfying (\ref{BDD PS}%
), where 
\begin{eqnarray}
&&
\begin{array}{l}
c:=\inf\limits_{\gamma \in \Gamma }\max\limits_{u\in \gamma \left( M\right)
}I\left( u\right) 
\end{array}
\nonumber \\
&&
\begin{array}{l}
\Gamma :=\{\gamma \in C\left( M,X\right) :\gamma _{\mid _{M_{0}}}=\mathrm{id}%
\},~~M:=\left\{ u\in Y:\left\| u\right\| \leq \rho \right\} .
\end{array}
\label{Thm(SD): def}
\end{eqnarray}
\end{thm}

\proof%
Let $M$ and $M_{0}$ be as in Theorem \ref{THM: SD} and set $\Gamma
_{0}:=\left\{ \gamma _{0}\right\} $ where $\gamma _{0}$ is the identity map
of $M_{0}$, so that $\Gamma $ of Theorem \ref{THM MAIN} reduces to the one
of (\ref{Thm(SD): def}). Let $\varepsilon _{0}\in \left( 0,\rho \right) $
and for every $u\in M$ define 
\[
\sigma \left( u\right) :=\left\{ 
\begin{array}{ll}
\dfrac{\rho }{\rho -\varepsilon _{0}}u\medskip & \text{if }\left\| u\right\|
\leq \rho -\varepsilon _{0} \\ 
\rho \dfrac{u}{\left\| u\right\| } & \text{if }\left\| u\right\| \geq \rho
-\varepsilon _{0}\,.
\end{array}
\right. 
\]
Since standard nonretractability arguments show that $\gamma \left( M\right)
\cap Z\neq \varnothing $ for every $\gamma \in \Gamma $, and thus $c\geq
\inf_{u\in Z}I\left( u\right) $, assumption (\ref{Thm(SD): hp}) imply (\ref
{a < c})-(\ref{sup max ->0}) and the conclusion then follows by applying
Theorem \ref{THM MAIN}.%
\endproof%
\bigskip

Denote $\mathbb{R}_{+}:=\left( 0,+\infty \right) $ and $\mathbb{R}_{+}^{0}:=\left[
0,+\infty \right) $.

\begin{thm}
\label{THM: L}Assume $X=Y\oplus Z$ with $\dim Y<\infty $ and assume that
there exist $\rho >r>0$ and $z\in Z$ with $\left\| z\right\| =r$ such that 
\[
\inf_{u\in Z,\left\| u\right\| =r}I\left( u\right) >\max_{u\in M_{0}}I\left(
u\right) \quad \text{and}\quad \lim\limits_{n\rightarrow \infty }\max_{u\in
M_{0}}\left\| \psi _{n}\left( u\right) -u\right\| =0
\]
where $M_{0}:=\left\{ u\in Y:\left\| u\right\| \leq \rho \right\} \cup
\left\{ u\in Y+\mathbb{R}_{+}z:\left\| u\right\| =\rho \right\} $. Then there
exists a sequence $\left\{ w_{n}\right\} \subset X$ satisfying (\ref{BDD PS}%
), where 
\begin{eqnarray*}
&&
\begin{array}{l}
c:=\inf\limits_{\gamma \in \Gamma }\max\limits_{u\in \gamma \left( M\right)
}I\left( u\right) 
\end{array}
\\
&&
\begin{array}{l}
\Gamma :=\{\gamma \in C\left( M,X\right) :\gamma _{\mid _{M_{0}}}=\mathrm{id}%
\},~~M:=\left\{ u\in Y+\mathbb{R}_{+}^{0}z:\left\| u\right\| \leq \rho \right\}
.
\end{array}
\end{eqnarray*}
\end{thm}

\proof%
One proceeds as for Theorem \ref{THM: SD} and the conclusion ensues from
applying Theorem \ref{THM MAIN}, provided that there exist $\varepsilon
_{0}>0$ and $\sigma \in C\left( M,M\right) $ satisfying \emph{i)}, \emph{ii)}%
, \emph{iii)}. In order to check this property, observe that $M=\left\{ u\in
Y+\mathbb{R}_{+}^{0}z:\left\| u\right\| \leq \rho \right\} $ is homeomorphic to
a finite dimensional compact ball and, for any fixed $\bar{\varepsilon}\in
\left( 0,1\right) $, define 
\[
\bar{\sigma}\left( p\right) :=\left\{ 
\begin{array}{ll}
\dfrac{1}{1-\bar{\varepsilon}}p\medskip & \text{if }\left| p\right| \leq 1-%
\bar{\varepsilon} \\ 
\dfrac{p}{\left| p\right| } & \text{if }\left| p\right| \geq 1-\bar{%
\varepsilon}
\end{array}
\right. \quad \forall p\in \bar{D}:=\left\{ p\in \mathbb{R}^{m+1}:\left|
p\right| \leq 1\right\} 
\]
where $m:=\dim Y+1$. Then it is easy to see that one can take $\sigma :=\phi
^{-1}\circ \bar{\sigma}\circ \phi $ where $\phi :Y\oplus \mathbb{R}z\rightarrow 
\mathbb{R}^{m+1}$ is any homeomorphism such that $\phi \left( M\right) =\bar{D}$
and $\varepsilon _{0}$ is any radius such that $\phi \left(
M_{0}^{\varepsilon _{0}}\right) \subseteq \left\{ p\in \bar{D}:\min_{q\in
\partial \bar{D}}\left| q-p\right| <\bar{\varepsilon}\right\} $.%
\endproof%

\section{Proof of Theorem \ref{THM MAIN}}

This section is devoted to the proof of Theorem \ref{THM MAIN}, which will
be achieved through several lemmas. Accordingly we assume all the hypotheses
of the theorem and, by (\ref{a < c}), fix any $\delta _{*}>0$ such that 
\begin{equation}
c_{0}:=\sup_{\gamma _{0}\in \Gamma _{0}}\max_{u\in \gamma _{0}\left(
M_{0}\right) }I\left( u\right) <c-\delta _{*}\,.  \label{a :=}
\end{equation}

\begin{lem}
\label{Lem: esiste eps_star}There exists $\varepsilon _{*}>0$ such that for
all $\gamma \in \Gamma $ and $u\in X$ one has 
\[
\mathrm{dist}\left( u,\gamma \left( M_{0}\right) \right) <\varepsilon
_{*}\,\Rightarrow \,I\left( u\right) <c_{0}+\delta _{*}\,.
\]
\end{lem}

\proof%
Let $C:=\left\{ u\in X:I\left( u\right) \geq c_{0}+\delta _{*}\right\} $, $%
K_{0}:=\bigcup_{\gamma _{0}\in \Gamma _{0}}\gamma _{0}\left( M_{0}\right) $
and $\varepsilon _{*}:=\mathrm{dist}\left( C,K_{0}\right) $ ($=\inf_{u\in
C,w\in K_{0}}\left\| u-w\right\| $). Recall that $K_{0}$ is compact by
assumption and assume by contradiction $\varepsilon _{*}=0$. Then there
exist $\left\{ u_{n}\right\} \subset C$ and $\left\{ w_{n}\right\} \subset
K_{0}$ such that $\left\| u_{n}-w_{n}\right\| \rightarrow 0$ and (up to a
subsequence) $w_{n}\rightarrow w_{0}\in K_{0}$, whence $\left\|
u_{n}-w_{0}\right\| \leq \left\| u_{n}-w_{n}\right\| +\left\|
w_{n}-w_{0}\right\| \rightarrow 0$ and thus $u_{n}\rightarrow w_{0}$. Since $%
C$ is closed, this gives $w_{0}\in C$, which is a contradiction because (\ref
{a :=}) implies $C\cap K_{0}=\varnothing $. Then the claim holds because $%
\mathrm{dist}\left( u,\gamma \left( M_{0}\right) \right) <\varepsilon _{*}$
and $I\left( u\right) \geq c_{0}+\delta _{*}$ would yield the existence of $%
w\in \gamma \left( M_{0}\right) \subseteq K_{0}$ (recall $\gamma _{\mid
_{M_{0}}}\in \Gamma _{0}$) such that $\left\| u-w\right\| <\varepsilon _{*}$
with $u\in C$.%
\endproof%
\bigskip

Henceforth, we fix any $A_{*}>0$ such that 
\[
\left| c\right| \liminf_{n\rightarrow \infty }\frac{\left| 1-\beta
_{n}\right| }{\alpha _{n}-\beta _{n}}<A_{*}<\infty . 
\]
Moreover, passing in case to subsequences of $\left\{ \psi _{n}\right\} $, $%
\left\{ \alpha _{n}\right\} $ and $\left\{ \beta _{n}\right\} $, by \textbf{(%
}$\mathbf{\Psi }$\textbf{{\tiny 2})} and (\ref{sup max ->0}) we assume that
for all $n$ there holds 
\begin{equation}
\alpha _{n}\left( \frac{\left| 1-\beta _{n}\right| }{\alpha _{n}-\beta _{n}}%
\left| c\right| +\left( 1+\beta _{n}\right) \left( \alpha _{n}-\beta
_{n}\right) \right) <A_{*}  \label{< A_0}
\end{equation}
and 
\begin{equation}
\sup_{\gamma _{0}\in \Gamma _{0}}\max_{u\in \gamma _{0}\left( M_{0}\right)
}\left\| \psi _{n}\left( u\right) -u\right\| <\frac{\varepsilon _{*}}{2}
\label{sup max <}
\end{equation}
where $\varepsilon _{*}>0$ is given by Lemma \ref{Lem: esiste eps_star}.

\begin{lem}
\label{Lem: I<a}For every $n$ and $\gamma \in \Gamma $ there exist $%
\varepsilon _{\gamma }>0$ and $\gamma ^{n}\in \Gamma $ such that for all $%
p\in M$ one has 
\begin{eqnarray}
&&
\begin{array}{l}
d\left( p,M_{0}\right) \geq \varepsilon _{\gamma }\,\Rightarrow \,\gamma
^{n}\left( p\right) =\psi _{n}\left( \gamma \left( \sigma \left( p\right)
\right) \right) 
\end{array}
\label{Lem( I<a ): >} \\
&&
\begin{array}{l}
d\left( p,M_{0}\right) <\varepsilon _{\gamma }\,\Rightarrow \,I\left( \gamma
^{n}\left( p\right) \right) <c_{0}+\delta _{*}\,.
\end{array}
\label{Lem( I<a ): <}
\end{eqnarray}
\end{lem}

\proof%
Let $\gamma \in \Gamma $. Since $\gamma $ is uniformly continuous on $M_{0}$
and $\sigma $ is uniformly continuous on $M_{0}^{\varepsilon _{0}}$ with $%
\sigma \left( M_{0}^{\varepsilon _{0}}\right) \subseteq M_{0}$, $\gamma
\circ \sigma $ turns out to be uniformly continuous on $M_{0}^{\varepsilon
_{0}}$ and thus there exists $\varepsilon _{\gamma }\in \left( 0,\varepsilon
_{0}\right) $ such that $p_{1},p_{2}\in M_{0}^{\varepsilon _{0}}$ and $%
d\left( p_{1},p_{2}\right) <\varepsilon _{\gamma }$ imply $\left\| \gamma
\left( \sigma \left( p_{1}\right) \right) -\gamma \left( \sigma \left(
p_{2}\right) \right) \right\| <\varepsilon _{*}/2$. For every $n$ and $p\in
M $ define 
\[
\gamma ^{n}\left( p\right) :=\left\{ 
\begin{array}{ll}
\psi _{n}(\gamma (\sigma \left( p\right) ))\medskip & \text{if }d\left(
p,M_{0}\right) \geq \varepsilon _{\gamma } \\ 
\mu \left( p\right) \psi _{n}(\gamma (\sigma \left( p\right) ))+\left( 1-\mu
\left( p\right) \right) \gamma \left( \sigma \left( p\right) \right) ~ & 
\text{if }d\left( p,M_{0}\right) \leq \varepsilon _{\gamma }
\end{array}
\right. 
\]
where $\mu :M\rightarrow \left[ 0,1\right] $ is any Uryson function such
that $\mu \left( p\right) =1$ if $d\left( p,M_{0}\right) =\varepsilon
_{\gamma }$ and $\mu \left( p\right) =0$ if $p\in M_{0}$. So $\gamma
^{n}:M\rightarrow X$ is continuous and satisfies $\gamma ^{n}\,_{\mid
_{M_{0}}}=\gamma \circ \sigma _{\mid _{M_{0}}}=\gamma _{\mid _{M_{0}}}\in
\Gamma _{0}$, namely, $\gamma ^{n}\in \Gamma $. Now let $p\in M$ satisfy $%
d\left( p,M_{0}\right) <\varepsilon _{\gamma }$ and let $\hat{p}\in M_{0}$
be such that $d\left( p,\hat{p}\right) <\varepsilon _{\gamma }$. Since $p,%
\hat{p}\in M_{0}^{\varepsilon _{0}}$, we deduce 
\[
\left\| \gamma \left( \sigma \left( p\right) \right) -\gamma \left( \hat{p}%
\right) \right\| =\left\| \gamma \left( \sigma \left( p\right) \right)
-\gamma \left( \sigma \left( \hat{p}\right) \right) \right\| <\frac{%
\varepsilon _{*}}{2} 
\]
which gives in turn 
\[
\left\| \psi _{n}\left( \gamma \left( \sigma \left( p\right) \right) \right)
-\gamma \left( \hat{p}\right) \right\| \leq \left\| \psi _{n}\left( \gamma
\left( \sigma \left( p\right) \right) \right) -\gamma \left( \sigma \left(
p\right) \right) \right\| +\frac{\varepsilon _{*}}{2}<\varepsilon _{*} 
\]
by (\ref{sup max <}). Therefore one infers $\left\| \gamma ^{n}\left(
p\right) -\gamma \left( \hat{p}\right) \right\| <\varepsilon _{*}$ by
convexity and then (\ref{Lem( I<a ): <}) follows from Lemma \ref{Lem: esiste
eps_star}.%
\endproof%

\begin{cor}
\label{Cor: c <= c_t}For all $n$ one has $c\leq \inf_{\gamma \in \Gamma
}\sup_{u\in \gamma \left( M\right) }I\left( \psi _{n}\left( u\right) \right) 
$.
\end{cor}

\proof%
Set $c_{n}:=\inf_{\gamma \in \Gamma }\sup_{u\in \gamma \left( M\right)
}I\left( \psi _{n}\left( u\right) \right) $ for brevity. Let $\delta >0$,
fix $\gamma \in \Gamma $ such that $\sup_{u\in \gamma \left( M\right)
}I\left( \psi _{n}\left( u\right) \right) \leq c_{n}+\delta $ and consider $%
\gamma ^{n}\in \Gamma $ given by Lemma \ref{Lem: I<a}. Then, by (\ref{Lem(
I<a ): <}), (\ref{a :=}) and (\ref{Lem( I<a ): >}), we have 
\begin{eqnarray*}
c &\leq &\sup_{u\in \gamma ^{n}\left( M\right) }I\left( u\right) =\sup_{p\in
M}I\left( \gamma ^{n}\left( p\right) \right) =\sup_{p\in M,d\left(
p,M_{0}\right) \geq \varepsilon _{\gamma }}I\left( \gamma ^{n}\left(
p\right) \right) \\
&=&\sup_{p\in M,d\left( p,M_{0}\right) \geq \varepsilon _{\gamma }}I\left(
\psi _{n}\left( \gamma \left( \sigma \left( p\right) \right) \right) \right)
\leq \sup_{p\in M}I\left( \psi _{n}\left( \gamma \left( \sigma \left(
p\right) \right) \right) \right) \\
&\leq &\sup_{p\in M}I\left( \psi _{n}\left( \gamma \left( p\right) \right)
\right) =\sup_{u\in \gamma \left( M\right) }I\left( \psi _{n}\left( u\right)
\right) \leq c_{n}+\delta
\end{eqnarray*}
which yields the claim as $\delta $ is arbitrary.%
\endproof%
\bigskip

From now on, to any $n$ we associate a map $\gamma _{n}\in \Gamma $
satisfying 
\begin{equation}
\sup_{u\in \gamma _{n}\left( M\right) }I\left( u\right) \leq c+\left( \alpha
_{n}-\beta _{n}\right) ^{2},  \label{max < c + beta}
\end{equation}
by which we define 
\[
\Lambda _{n}:=\left\{ u\in \gamma _{n}\left( M\right) :I\left( \psi
_{n}\left( u\right) \right) \geq c-\left( \alpha _{n}-\beta _{n}\right)
^{2}\right\} . 
\]
Note that $\Lambda _{n}\neq \varnothing $ by Corollary \ref{Cor: c <= c_t},
which implies $\sup_{u\in \gamma _{n}\left( M\right) }I\left( \psi
_{n}\left( u\right) \right) >c-\left( \alpha _{n}-\beta _{n}\right) ^{2}$
(indeed, by continuity, $\Lambda _{n}$ even contains an open subset of $%
\gamma _{n}\left( M\right) $).

\begin{lem}
\label{Lem: A(u) <}For every $n$ and $u\in \Lambda _{n}$ one has $A\left(
u\right) \leq A_{*}$.
\end{lem}

\proof%
Fix any $u\in \Lambda _{n}$. From \textbf{(}$\mathbf{\Psi }$\textbf{{\tiny 1}%
)} it follows that 
\begin{eqnarray*}
I\left( \psi _{n}\left( u\right) \right) &=&A\left( \psi _{n}\left( u\right)
\right) -B\left( \psi _{n}\left( u\right) \right) \leq \frac{1}{\alpha _{n}}%
A\left( u\right) -\frac{1}{\beta _{n}}B\left( u\right) \\
&=&\frac{1}{\alpha _{n}}A\left( u\right) -\frac{1}{\beta _{n}}\left( A\left(
u\right) -I\left( u\right) \right) =-\frac{\alpha _{n}-\beta _{n}}{\alpha
_{n}\beta _{n}}A\left( u\right) +\frac{1}{\beta _{n}}I\left( u\right)
\end{eqnarray*}
which yields, by (\ref{max < c + beta}) and the definition of $\Lambda _{n}$%
, 
\[
\frac{\alpha _{n}-\beta _{n}}{\alpha _{n}\beta _{n}}A\left( u\right) \leq 
\frac{1}{\beta _{n}}I\left( u\right) -I\left( \psi _{n}\left( u\right)
\right) \leq \frac{1}{\beta _{n}}\left( c+\left( \alpha _{n}-\beta
_{n}\right) ^{2}\right) -\left( c-\left( \alpha _{n}-\beta _{n}\right)
^{2}\right) 
\]
and thus 
\begin{eqnarray*}
A\left( u\right) &\leq &\alpha _{n}\left( \frac{1-\beta _{n}}{\alpha
_{n}-\beta _{n}}c+\left( 1+\beta _{n}\right) \left( \alpha _{n}-\beta
_{n}\right) \right) \\
&\leq &\alpha _{n}\left( \frac{\left| 1-\beta _{n}\right| }{\alpha
_{n}-\beta _{n}}\left| c\right| +\left( 1+\beta _{n}\right) \left( \alpha
_{n}-\beta _{n}\right) \right) <A_{*}
\end{eqnarray*}
by (\ref{< A_0}).%
\endproof%

\begin{rem}
By suitably changing assumption (\ref{< A_0}), the above proof of Lemma \ref
{Lem: A(u) <} shows that neither that case $c>0$ and $\liminf_{n\rightarrow
\infty }\left( 1-\beta _{n}\right) /\left( \alpha _{n}-\beta _{n}\right) <0$%
, nor the case $c<0$ and $\limsup_{n\rightarrow \infty }\left( 1-\beta
_{n}\right) /\left( \alpha _{n}-\beta _{n}\right) >0$ can occur under the
assumed hypotheses, since this would yield the contradiction $A\left(
u\right) <0$ for all $u\in \Lambda _{n}\,\left( \neq \varnothing \right) $
with $n$ large enough.
\end{rem}

\begin{lem}
\label{Lem: cece}One has $\lim\limits_{n\rightarrow \infty
}\sup\limits_{u\in \Lambda _{n}}\left| I\left( u\right) -c\,\right|
=\lim\limits_{n\rightarrow \infty }\sup\limits_{u\in \Lambda _{n}}\left|
I\left( \psi _{n}\left( u\right) \right) -I\left( u\right) \right| =0$.
\end{lem}

\proof%
Since $A\geq 0$, for every $u\in \Lambda _{n}$ from \textbf{(}$\mathbf{\Psi }
$\textbf{{\tiny 1})} one deduces 
\[
I\left( \psi _{n}\left( u\right) \right) \leq \frac{1}{\alpha _{n}}A\left(
u\right) -\frac{1}{\beta _{n}}B\left( u\right) =-\frac{\alpha _{n}-\beta _{n}%
}{\alpha _{n}\beta _{n}}A\left( u\right) +\frac{1}{\beta _{n}}I\left(
u\right) \leq \frac{1}{\beta _{n}}I\left( u\right) 
\]
so that, by the definition of $\Lambda _{n}$, one has 
\[
I\left( u\right) -c\geq \beta _{n}I\left( \psi _{n}\left( u\right) \right)
-c\geq \beta _{n}\left( c-\left( \alpha _{n}-\beta _{n}\right) ^{2}\right)
-c\,. 
\]
By (\ref{max < c + beta}) this yields $\sup_{u\in \Lambda _{n}}\left|
I\left( u\right) -c\right| \rightarrow 0$ as $n\rightarrow \infty $. By (\ref
{max < c + beta}) and the definition of $\Lambda _{n}$ one also has $I\left(
u\right) -I\left( \psi _{n}\left( u\right) \right) \leq 2\left( \alpha
_{n}-\beta _{n}\right) ^{2}$ for all $u\in \Lambda _{n}$, whereas \textbf{(}$%
\mathbf{\Psi }$\textbf{{\tiny 1})} and Lemma \ref{Lem: A(u) <} yield 
\begin{eqnarray*}
I\left( \psi _{n}\left( u\right) \right) -I\left( u\right) &=&A\left( \psi
_{n}\left( u\right) \right) -A\left( u\right) +B\left( \psi _{n}\left(
u\right) \right) +B\left( u\right) \\
&\leq &\frac{\left| 1-\alpha _{n}\right| }{\alpha _{n}}A\left( u\right) +%
\frac{\left| 1-\beta _{n}\right| }{\beta _{n}}\left| B\left( u\right) \right|
\\
&\leq &\frac{\left| 1-\alpha _{n}\right| }{\alpha _{n}}A_{*}+\frac{\left|
1-\beta _{n}\right| }{\beta _{n}}\left( A\left( u\right) +\left| I\left(
u\right) \right| \right) \\
&\leq &\left( \frac{\left| 1-\alpha _{n}\right| }{\alpha _{n}}+\frac{\left|
1-\beta _{n}\right| }{\beta _{n}}\right) A_{*}+\frac{\left| 1-\beta
_{n}\right| }{\beta _{n}}\left( c+\sup_{u\in \Lambda _{n}}\left| I\left(
u\right) -c\right| \right) .
\end{eqnarray*}
Hence $\sup_{u\in \Lambda _{n}}\left| I\left( \psi _{n}\left( u\right)
\right) -I\left( u\right) \right| \rightarrow 0$ as $n\rightarrow \infty $.%
\endproof%
\bigskip

Now fix $l_{1},l_{2}>0$ such that $A_{*}<l_{1}<l_{2}$. For every $k\geq 1$,
define 
\[
U_{k}:=\left\{ u\in X:A\left( u\right) \leq l_{2}+\frac{1}{k},~\left|
I\left( u\right) -c\right| \leq \frac{1}{k}\right\} 
\]
and, by \textbf{(}$\mathbf{\Psi }$\textbf{{\tiny 2})} and Lemma \ref{Lem:
cece}, take $n_{k}\in \mathbb{N}$ such that 
\begin{equation}
\alpha _{n_{k}}>\frac{A_{*}}{l_{1}}\,,\quad \sup_{u\in \Lambda
_{n_{k}}}I\left( \psi _{n_{k}}\left( u\right) \right) \leq c+\frac{1}{16k}%
\,,\quad \sup_{u\in \Lambda _{n_{k}}}\left| I\left( u\right) -c\right| \leq 
\frac{1}{k}.  \label{n_k}
\end{equation}
Hence $\Lambda _{n_{k}}\subseteq U_{k}$ for all $k$ and $U_{k}$ is not empty.

In order to conclude the proof of Theorem \ref{THM MAIN} we need to apply a
well known deformation lemma (see \cite[Lemma 2.3]{Willem}), which, for
completeness, we recall here for the space $X$, our functional $I$ and its
minimax level $c$.

\begin{lem}
\label{Lem: deformation}Let $\mathcal{S}\subset X$ and $\varepsilon ,\delta
>0$ be such that $\left\| I^{\prime }\left( u\right) \right\| _{X^{\prime
}}\geq 8\varepsilon /\delta $ for all $u\in \mathcal{S}_{2\delta }$
satisfying $\left| I\left( u\right) -c\right| \leq 2\varepsilon $, where $%
\mathcal{S}_{2\delta }:=\left\{ v\in X:\mathrm{dist}\left( v,\mathcal{S}%
\right) \leq 2\delta \right\} $. Then there exists $\eta \in C\left( \left[
0,1\right] \times X,X\right) $ such that

\begin{itemize}
\item  $\eta \left( \tau ,\cdot \right) $ is a homeomorphism of $X$ for
every $\tau \in \left[ 0,1\right] $

\item  $\eta \left( \tau ,u\right) =u$ provided that $\tau =0$ or $\left|
I\left( u\right) -c\right| >2\varepsilon $ or $u\notin \mathcal{S}_{2\delta }
$

\item  $I\left( \eta \left( 1,u\right) \right) \leq c-\varepsilon $ provided
that $I\left( u\right) \leq c+\varepsilon $ and $u\in \mathcal{S}$

\item  $I\left( \eta \left( \cdot ,u\right) \right) $ is nonincreasing for
every $u\in X\,.\medskip $
\end{itemize}
\end{lem}

\proof[Proof of Theorem \ref{THM MAIN}]
Consider the sublevels $A^{l_{1}},A^{l_{2}}$ of $A$ and, by hypothesis \textbf{(}$\mathbf{%
A}$\textbf{)}, fix $\delta _{0}>0$ such that 
\begin{equation}
\forall u\in X\qquad \mathrm{dist}\left( u,A^{l_{1}}\right) \leq \delta
_{0}~\Rightarrow ~u\in A^{l_{2}}\,.  \label{semicont}
\end{equation}
Then assume by contradiction that 
\begin{equation}
\exists \bar{k}>\max \left\{ \frac{1}{\delta _{0}^{2}},\frac{1}{8\left(
c-c_{0}-\delta _{*}\right) }\right\} \quad \forall u\in U_{\bar{k}}\quad
\left\| I^{\prime }\left( u\right) \right\| _{X^{\prime }}\geq \frac{1}{%
\sqrt{\bar{k}}}\,.  \label{hp assurda}
\end{equation}
Recall that $c-c_{0}-\delta _{*}>0$ by (\ref{a :=}). We apply Lemma \ref
{Lem: deformation} with $\mathcal{S}=A^{l_{1}}$, $\varepsilon =1/16\bar{k}$
and $\delta =1/2\sqrt{\bar{k}}$ (so that $8\varepsilon /\delta =1/\sqrt{\bar{%
k}}$). Observe that if $u\in \mathcal{S}_{2\delta }$ satisfies $\left|
I\left( u\right) -c\right| \leq 1/8\bar{k}$ then $u\in U_{\bar{k}}$ (and
thus the last inequality of (\ref{hp assurda}) holds), since $2\delta =1/%
\sqrt{\bar{k}}<\delta _{0}$ and (\ref{semicont}) applies. So there exists a
homeomorphism $\Phi :X\rightarrow X$ (namely $\Phi :=\eta \left( 1,\cdot
\right) $ of Lemma \ref{Lem: deformation}) such that

\begin{itemize}
\item[(i)]  $\Phi \left( u\right) =u$ if $\left| I\left( u\right) -c\right|
\geq c-c_{0}-\delta _{*}$ (recall $c-c_{0}-\delta _{*}>\frac{1}{8\bar{k}}%
=2\varepsilon $)

\item[(ii)]  $I\left( \Phi \left( u\right) \right) \leq c-\frac{1}{16\bar{k}}
$ if $A\left( u\right) \leq l_{1}$ and $I\left( u\right) \leq c+\frac{1}{16%
\bar{k}}$

\item[(iii)]  $I\left( \Phi \left( u\right) \right) \leq I\left( u\right) $
for every $u\in X\,$,
\end{itemize}

\noindent by which we define the mapping $\gamma :=\Phi \circ \gamma _{n_{%
\bar{k}}}^{n_{\bar{k}}}\in C\left( M,X\right) $, where $\gamma _{n_{\bar{k}%
}}^{n_{\bar{k}}}=(\gamma _{n_{\bar{k}}})^{n_{\bar{k}}}\in \Gamma $ is the
mapping associated to $\gamma _{n_{\bar{k}}}$ by Lemma \ref{Lem: I<a}. Then,
by (\ref{Lem( I<a ): <}) of Lemma \ref{Lem: I<a}, $p\in M$ and $d\left(
p,M_{0}\right) <\varepsilon _{\gamma _{n_{\bar{k}}}}$ imply $I\left( \gamma
_{n_{\bar{k}}}^{n_{\bar{k}}}\left( p\right) \right) <c_{0}+\delta _{*}$,
which yields 
\[
\left| \,I\left( \gamma _{n_{\bar{k}}}^{n_{\bar{k}}}\left( p\right) \right)
-c\,\right| =c-I\left( \gamma _{n_{\bar{k}}}^{n_{\bar{k}}}\left( p\right)
\right) >c-c_{0}-\delta _{*} 
\]
and so, by (i), $\gamma \left( p\right) =\Phi \left( \gamma _{n_{\bar{k}%
}}^{n_{\bar{k}}}\left( p\right) \right) =\gamma _{n_{\bar{k}}}^{n_{\bar{k}%
}}\left( p\right) $ and $I\left( \Phi \left( \gamma _{n_{\bar{k}}}^{n_{\bar{k%
}}}\left( p\right) \right) \right) =I\left( \gamma _{n_{\bar{k}}}^{n_{\bar{k}%
}}\left( p\right) \right) <c_{0}+\delta _{*}<c$. Hence, in particular, $%
\gamma _{\mid _{M_{0}}}=\gamma _{n_{\bar{k}}}^{n_{\bar{k}}}\,_{\mid
_{M_{0}}}\in \Gamma _{0}$ and thus $\gamma \in \Gamma $. Moreover, by (\ref
{Lem( I<a ): >}), one has 
\begin{eqnarray}
\sup_{u\in \gamma \left( M\right) }I\left( u\right) &=&\sup_{p\in M}I\left(
\Phi \left( \gamma _{n_{\bar{k}}}^{n_{\bar{k}}}\left( p\right) \right)
\right) =\sup_{p\in M,\,d\left( p,M_{0}\right) \geq \varepsilon _{\gamma
_{n_{\bar{k}}}}}I\left( \Phi \left( \gamma _{n_{\bar{k}}}^{n_{\bar{k}%
}}\left( p\right) \right) \right)  \nonumber \\
&=&\sup_{p\in M,\,d\left( p,M_{0}\right) \geq \varepsilon _{\gamma _{n_{\bar{%
k}}}}}I\left( \Phi \left( \psi _{n_{\bar{k}}}\left( \gamma _{n_{\bar{k}%
}}\left( \sigma \left( p\right) \right) \right) \right) \right)  \nonumber \\
&\leq &\sup_{p\in M}I\left( \Phi \left( \psi _{n_{\bar{k}}}\left( \gamma
_{n_{\bar{k}}}\left( \sigma \left( p\right) \right) \right) \right) \right) 
\nonumber \\
&\leq &\sup_{u\in \gamma _{n_{\bar{k}}}\left( M\right) }I\left( \Phi \left(
\psi _{n_{\bar{k}}}\left( u\right) \right) \right)  \label{sup <}
\end{eqnarray}
and this yields a contradiction. Indeed, on one hand, if $u\in \gamma _{n_{%
\bar{k}}}\left( M\right) \setminus \Lambda _{n_{\bar{k}}}$, by (iii) and the
definition of $\Lambda _{t_{\bar{m}}}$ we get 
\[
I\left( \Phi \left( \psi _{n_{\bar{k}}}\left( u\right) \right) \right) \leq
I\left( \psi _{n_{\bar{k}}}\left( u\right) \right) <c-\left( \alpha _{n_{%
\bar{k}}}-\beta _{n_{\bar{k}}}\right) . 
\]
On the other hand, if $u\in \Lambda _{n_{\bar{k}}}$, we have $I(\psi _{n_{%
\bar{k}}}\left( u\right) )\leq c+1/16\bar{k}$ by (\ref{n_k}) and $A(\psi
_{n_{\bar{k}}}\left( u\right) )\leq A\left( u\right) /\alpha _{n_{\bar{k}%
}}<l_{1}$ by \textbf{(}$\mathbf{\Psi }$\textbf{{\tiny 1})}, (\ref{n_k}) and
Lemma \ref{Lem: A(u) <}, so that (ii) gives 
\[
I\left( \Phi \left( \psi _{n_{\bar{k}}}\left( u\right) \right) \right) \leq
c-\frac{1}{16\bar{k}}\,. 
\]
Therefore, by (\ref{sup <}), one obtains $\sup_{u\in \gamma \left( M\right)
}I\left( u\right) <c$ and the definition (\ref{a < c}) of $c$ is
contradicted.%
\endproof%


\begin{thebibliography}{99}
\bibitem{Ambr-Rab}  \textsc{Ambrosetti A., Rabinowitz P.H.}, \emph{Dual
variational methods in critical point theory and applications}, J. Funct.
Anal. \textbf{14} (1973), 349-381.

\bibitem{Ambr-Ruiz}  \textsc{Ambrosetti A., Ruiz D.}, \emph{Multiple bound
states for the Schr\"{o}dinger-Poisson problems}, Comm. Contemp. Math. 
\textbf{10} (2008), 1-14.

\bibitem{Ambr-Str}  \textsc{Ambrosetti A., Struwe M.}, \emph{Existence of
steady vortex rings in an ideal fluid}, Arch. Rational Mech. Anal. \textbf{%
108} (1989), 97-109.

\bibitem{BBR 3}  \textsc{Badiale M., Benci V., Rolando S.}, \emph{Three
dimensional vortices in the nonlinear wave equation}, Boll. Unione Mat.
Ital., Serie IX, \textbf{2} (2009), 105-134.

\bibitem{BR TMA}  \textsc{Badiale M., Rolando S.}, \emph{Nonlinear elliptic
equations with subhomogeneous potentials}, Nonlinear Anal. \textbf{72}
(2010), 602-617.

\bibitem{GR curves}  \textsc{Guida M., Rolando S.}, \emph{Symmetric }$\kappa 
$\emph{-loops}, Diff. Int. Equations \textbf{23} (2010), 861-898.

\bibitem{GR bdd}  \textsc{Guida M., Rolando S.}, \emph{On the existence of
bounded Palais-Smale sequences and applications to nonlinear equations
without superlinearity assumptions}, work in progress.

\bibitem{Jabri}  \textsc{Jabri Y.}, \emph{The Mountain Pass Theorem},
Cambridge University Press, 2003.

\bibitem{JJ}  \textsc{Jeanjean L.}, \emph{On the existence of bounded
Palais-Smale sequences and application to a Landesman-Lazer type problem set
on }$\mathbb{R}^{N}$, Proc. Roy. Soc. Edinburgh A \textbf{129} (1999), 787-809.

\bibitem{JJ 2}  \textsc{Jeanjean L.}, \emph{Local condition insuring
bifurcation from the continuous spectrum}, Math. Z. \textbf{232} (1999),
651-664.

\bibitem{JJ-Tanaka}  \textsc{Jeanjean L., Tanaka K.}, \emph{A positive
solution for a nonlinear Schr\"{o}dinger equation on }$\mathbb{R}^{N}$, Indiana
Univ. Math. J. \textbf{54} (2005), 443-464.

\bibitem{JJ-Toland}  \textsc{Jeanjean L., Toland J.F.}, \emph{Bounded
Palais-Smale Mountain-Pass sequences}, C.R. Acad. Sci. Paris, Ser. I \textbf{%
327} (1998), 23-28.

\bibitem{Miya-Souto}  \textsc{Miyagaki O.H., Souto M.A.S.}, \emph{%
Superlinear problems without Ambrosetti and Rabinowitz growth condition}, J.
Differential Equations \textbf{245} (2008), 3628-3638.

\bibitem{Rabier}  \textsc{Rabier P.J.}, \emph{Bounded Palais-Smale sequences
for functionals with a mountain pass geometry}, Arch. Math. \textbf{88}
(2007), 143-152.

\bibitem{Schech-Zou}  \textsc{Schechter M., Zou W.}, \emph{Weak linking
theorems and Schr\"{o}dinger equations with critical Sobolev exponent},
Control, Optimisation and Calculus of Variations \textbf{9} (2003), 601-619.

\bibitem{Schech-Zou 04}  \textsc{Schechter M., Zou W.}, \emph{Superlinear
problems}, Pacific J. Math. \textbf{214} (2004), 145-160.

\bibitem{Strick}  \textsc{Struwe M.}, \emph{The existence of surfaces of
constant mean curvature with free boundaries}, Acta Math. \textbf{160}
(1988), 19-64.

\bibitem{Szulkin-Zou}  \textsc{Szulkin A., Zou W.}, \emph{Homoclinic orbits
for asymptotically linear hamiltonian systems}, J. Funct. Anal. \textbf{187}
(2001), 25-41.

\bibitem{Wang-Tang}  \textsc{Wang J., Tang C.L.}, \emph{Existence and
multiplicity of solutions for a class of superlinear }$p$\emph{-laplacian
equations}, Boundary Value Problems, vol. 2006, Article ID 47275, 12 pages,
2006. doi:10.1155/BVP/2006

\bibitem{Willem}  \textsc{Willem M.}, \emph{Minimax Theorems}, PNLDE, vol.
24, Birkh\"{a}user, Boston 1996.

\bibitem{Will-Zou}  \textsc{Willem M., Zou W.}, \emph{On a Schr\"{o}dinger
equation with periodic potential and spectrum point zero}, Indiana Univ.
Math. J. \textbf{52} (2003), 109-132.

\bibitem{Zhou 98}  \textsc{Zhou H.S.}, \emph{Positive solution for a
semilinear elliptic equation which is almost linear at infinity}, Z. Angew.
Math. Phys. \textbf{49} (1998), 896-906.

\bibitem{Zhou}  \textsc{Zhou H.S.}, \emph{Existence of asymptotically linear
Dirichlet problem}, Nonlinear Analysis \textbf{44} (2001), 909-918.

\bibitem{Zou}  \textsc{Zou W.}, \emph{Variant fountain theorems and their
applications}, Indiana Univ. Math. J. \textbf{52} (2003), 109-132.
\end{thebibliography}
\end{document}